\documentclass[10pt,a4paper]{article}
\usepackage[margin=1in]{geometry}
\usepackage{amsmath,amsthm,amssymb,enumerate, bbm ,graphicx,color,caption,upgreek, float, tikz, wasysym, subcaption,booktabs,longtable, appendix,graphics, pdfpages,rotating,mathtools,textcomp, array, bm}
\DeclarePairedDelimiter\floor{\lfloor}{\rfloor}
\usepackage[pdftex]{hyperref} 
\usepackage[customcolors,shade]{hf-tikz}

\usetikzlibrary{arrows, positioning,decorations.pathreplacing}
\definecolor{blauw}{RGB}{61,158,255}
\definecolor{donkerblauw}{RGB}{0,0,255}
\definecolor{donkergroen}{RGB}{46,148,0}
\definecolor{donkerrood}{RGB}{204,0,0}

\makeatletter 
\newcommand\mynobreakpar{\par\nobreak\@afterheading} 
\makeatother

\newcommand{\N}{\mathbb{N}}
\newcommand{\Z}{\mathbb{Z}}
\newcommand{\C}{\mathbb{C}}
\newcommand{\R}{\mathbb{R}}

\newcommand{\F}{\mathbb{F}}
\usepackage[english]{babel}

\newtheorem{theorem}{Theorem}[section]

\newtheorem{claim}[theorem]{Claim}
\newtheorem{proposition}[theorem]{Proposition}

\theoremstyle{definition}

\newtheorem*{examp*}{Example}

\makeatletter
\let\@fnsymbol\@arabic
\makeatother

\theoremstyle{plain}

\newfloat{Algorithm}{!hbt}{alg}

\newcommand{\dov}{\textbf{do }}
\newcommand{\ndv}{\textbf{end }}
\newcommand{\ifv}{\textbf{if }}

\newcommand{\printv}{\textbf{print }}

\newcommand{\foreachv}{\textbf{foreach }}

\newcounter{thm}[section]

\title{Semidefinite programming bounds for constant weight codes}
\date{}
\author{Sven Polak\thanks{Korteweg-De Vries Institute for Mathematics, University of Amsterdam. E-mail: \href{mailto:s.c.polak@uva.nl}{s.c.polak@uva.nl}. The research leading to these
results has received funding from the European Research Council under the European Union’s Seventh Framework Programme (FP7/2007-2013) / ERC grant agreement \textnumero 339109.}}
\begin{document}
\maketitle
\setcounter{footnote}{1}

\noindent \textbf{Abstract.} For nonnegative integers~$n,d,w$, let~$A(n,d,w)$ be the maximum size of a code~$C \subseteq \mathbb{F}_2^n$ with constant weight~$w$ and minimum distance at least~$d$. We consider two semidefinite programs based on quadruples of code words that yield several new upper bounds on~$A(n,d,w)$. The new upper bounds imply that~$A(22,8,10)=616$ and~$A(22,8,11)=672$. Lower bounds on~$A(22,8,10)$ and~$A(22,8,11)$ are obtained from the~$(n,d)=(22,7)$ shortened Golay code of size~$2048$. It can be concluded that the shortened Golay code is a union of constant weight~$w$ codes of sizes~$A(22,8,w)$. 

\,$\phantom{0}$

\noindent {\bf Keywords:} constant-weight code, upper bounds, semidefinite programming, Delsarte, Golay.

\noindent {\bf MSC 2010:} 94B65, 05E10, 90C22, 20C30.
\section{Introduction}

\noindent Let~$\F_2:=\{0,1\}$ denote the field of two elements and fix~$n \in \N$. A \emph{word} is an element~$v \in \F_2^n$. For two words~$u,v \in \F_2^n$, their \emph{(Hamming) distance}~$d_H(u,v)$ is the number of~$i$ with~$u_i \neq v_i$. A \emph{code} is a subset of~$\F_2^n$. For any code~$C \subseteq \F_2^n$, the \emph{minimum distance}~$d_{\text{min}}(C)$ of~$C$ is the minimum distance between any pair of distinct code words in~$C$. The \emph{weight} $\text{wt}(v)$ of a word~$v \in \F_2^n$ is the number of nonzero entries of~$v$. This paper considers \emph{constant weight codes}, i.e., codes in which all code words have a fixed weight~$w$. Then $A(n,d,w)$ is defined as the maximum size of a code of minimum distance at least~$d$ in which every codeword has weight~$w$. Moreover,~$A(n,d)$ is the maximum size of a code of minimum distance at least~$d$. Determining~$A(n,d)$ and~$A(n,d,w)$ for~$n,d,w\in \N$ are long-time focuses in combinatorial coding theory (cf.~$\cite{sloane}$).

In this paper we will consider two semidefinite programming upper bounds on~$A(n,d,w)$. Both upper bounds sharpen the classical Delsarte linear programming bound~$\cite{delsarte}$, as well as Schrijver's semidefinite programming bound based on a block diagonalization of the Terwilliger algebra~$\cite{schrijver}$.

The paper serves the following purposes. Firstly, the bound~$A_k(n,d)$ for unrestricted (non-constant weight) codes from~$\cite{semidef}$ is adapted to a bound~$A_k(n,d,w)$ for constant weight codes. Subsequently, a relaxation~$B_k(n,d,w)$ is formulated, which might also be of interest for unrestricted binary codes. By studying~$A_4(n,d,w)$ and~$B_4(n,d,w)$,  a sharpening of the Schrijver bound~$\cite{schrijver}$ for constant weight codes is obtained that is in most cases sharper than the bound from~$\cite{KT}$ (in which linear inequalities were added to the Schrijver bound). The constructed semidefinite programs are very large, but a symmetry reduction (using representation theory of the symmetric group) is given to reduce them to polynomial size. This finally leads to many new upper bounds on~$A(n,d,w)$, including the exact values~$A(22,8,10)=616$ and~$A(22,8,11)=672$.  

The once shortened Golay code, which is an~$(n,d)=(22,7)$-code of size 2048, contains the following numbers of words of a given weight~$w$ (and no words of other weights).
\begin{table}[H]
\centering
\begin{tabular}{l|lllllll}
weight~$w$ & 0 & 7 & 8 & 11 &12 &15 &16 \\\hline
$\# $ words  &1  & 176 & 330 & 672 & 616 &176 &77
\end{tabular}
\caption{Number of words of a given weight~$w$ contained in the shortened Golay code.}
\label{UB}
\end{table}
\noindent While it was already known that~$A(22,8,7)=176$,~$A(22,8,8)=330$ and~$A(22,8,6)=77$  (here note that~$A(n,d,w)=A(n,d,n-w)$), the results of this paper imply that~$A(22,8,10)=616$ and~$A(22,8,11)=672$. So if one collects all words of a given weight~$w$ in the shortened Golay code, the resulting code is a constant weight code of maximum size. The shortened Golay code therefore is a union of constant weight~$w$ codes of sizes~$A(22,8,w)$. The value~$A(22,8,10)=616$ together with the already known values implies that also the twice shortened extended Golay code (which is an~$(n,d)=(22,8)$-code of size~1024) has this property, since it contains~$1$,~$330$,~$616$ and~$77$ words of weight~$0$,~$8$,~$12$ and~$16$ (respectively).   It was already known that the Golay code, the extended Golay code and the once shortened extended Golay code have this property, i.e., that they are unions of constant weight codes of sizes~$A(n,d,w)$.

 Many tables with bounds on~$A(n,d,w)$ have been given in the literature~$\cite{agrell, table2, table3}$.  Tables with best currently known upper and lower bounds can be found on the website of Andries Brouwer~$\cite{brouwertable}$.

\subsection{The upper bounds~\texorpdfstring{$A_k(n,d,w)$}{A4(n,d,w)} and~\texorpdfstring{$B_k(n,d,w)$}{B4(n,d,w)}.}

We describe two upper bounds on~$A(n,d,w)$ based on quadruples of code words.  Fix~$n,d,w \in \N$ and let~$F \subseteq \mathbb{F}_2^n$ be the set of all words of constant weight~$w$. For~$k \in \Z_{\geq 0}$, let~$\mathcal{C}_k$ be the collection of codes~$C \subseteq F$ with~$|C|\leq k$.  For any  natural number~$j\leq k$ and~$D \in \mathcal{C}_j$, we define
\begin{align} 
\mathcal{C}_j(D) := \{C \in \mathcal{C}_j \,\, | \,\, C \supseteq D, \, |D|+2|C\setminus D| \leq j  \}.  
\end{align} 
Note that then~$|C \cup C'| \leq j$ for all~$C,C' \subseteq \mathcal{C}_j(D)$.   
Furthermore, for any function~$x : \mathcal{C}_k \to \R$ and~$D \in \mathcal{C}_j$ we define the~$\mathcal{C}_j(D) \times \mathcal{C}_j(D)$-matrix~$M_{j,D}(x)$ by $M_{j,D}(x)_{C,C'} : = x(C \cup C')$. Define the following number, which is an adaptation to constant weight codes of the bound~$A_k(n,d)$ from Gijswijt,  Mittelmann and Schrijver~\cite{semidef}:
\begin{align} \label{A4ndw}
A_k(n,d,w):=    \max \{ \sum_{v \in F} x(\{v\})\,\, &|\,\,x:\mathcal{C}_k \to \R, \,\, x(\emptyset )=1, x(S)=0 \text{ if~$d_{\text{min}}(S)<d$},\\[-1.1em]
& \,\, M_{k,D}(x) \text{ is positive semidefinite for each~$D$ in~$\mathcal{C}_k$}\}. \notag 
\end{align}
In this paper, we first consider~$A_4(n,d,w)$. Even after reductions (see the next subsection), the semidefinite program for computing~$A_{4}(n,d,w)$ is large in practice, although~$A_k(n,d,w)$ can be computed in polynomial time for fixed~$k$. In computing~$A_k(n,d,w)$, the matrix blocks coming from the matrices~$M_{k,D}(x)$ for~$D = \emptyset$  if~$k$ is even, and for~$|D|=1$ if~$k$ is odd, are often larger (in size and importantly, more variables occur in each matrix entry, yielding large semidefinite programs) than the blocks coming from~$M_{k,D}(x)$ for~$D$ of other cardinalities. This observation gives rise to the following relaxation of~$A_k(n,d,w)$, which is sharper than~$A_{k-1}(n,d,w)$ for~$k\geq 4$ (while it equals~$A_{k-1}(n,d,w)$ for~$k=3$, so assume~$k \geq 3$ in the following definition).
\begin{align} \label{Bndw}
B_k(n,d,w):=   \max \{ \sum_{v \in F} x(\{v\})\,\, &|\,\,x:\mathcal{C}_k \to \R, \,\, x(\emptyset )=1, x(S)=0 \text{ if~$d_{\text{min}}(S)<d$},\\[-1.1em]
&\,\, M_{k-1,D}(x) \text{ is positive semidefinite for each~$D \in \mathcal{C}_{k-1}$ with~$|D|<2$}, \notag 
\\
&\,\, M_{k,D}(x) \text{ is positive semidefinite for each~$D \in \mathcal{C}_{k}$ with~$|D|\geq 2$}\}. \notag 
\end{align}

\begin{proposition}
Fix~$k \geq 3$. For all~$n,d,w \in \N$, we have~$ A_{k-1}(n,d,w) \geq B_k(n,d,w) \geq A_k(n,d,w) \geq A(n,d)$.
\end{proposition}
\proof
Let~$C \subseteq \mathbb{F}_2^n$ be a constant weight~$w$ code with~$d_{\text{min}}(C)\geq d$ and~$|C| = A(n,d,w)$. Define~$x: \mathcal{C}_k \to \R$ by~$x(S)=1$ if~$S \subseteq C$ and~$x(S)=0$ else. Then~$x$ satisfies the conditions in~$(\ref{A4ndw})$, where the last condition is satisfied since~$M_{k,D}(x)_{C,C'}=x(C)x(C')$ for all~$C,C'\in \mathcal{C}_k$. Moreover, the objective value equals~$\sum_{v \in F} x(\{v\}) =|C|=A(n,d,w)$, which gives~$A_k(n,d,w) \geq A(n,d,w)$.

It is not hard to see that~$B_k(n,d,w) \geq A_k(n,d,w)$, as all constraints in~$(\ref{Bndw})$ follow from~$(\ref{A4ndw})$. Similarly, it follows that~$A_{k-1}(n,d,w)\geq B_k(n,d,w)$. 
\endproof

\noindent This paper considers~$A_4(n,d,w)$ and~$B_4(n,d,w)$, that is,~$k=4$. By symmetry we assume throughout that~$w \leq n/2$ (otherwise, add the all-ones word to each word in~$\mathbb{F}_2^n$ and replace~$w$ by~$n-w$). For computing~$A_4(n,d,w)$, it suffices to require that the matrices~$M_{k,D}(x)$ with~$|D|$ even are positive semidefinite. To see this, note that if~$D \subseteq C$ with~$|D|$ even and~$|C|=|D|+1$, then~$\mathcal{C}_k(C) \subseteq \mathcal{C}_k(D)$, i.e.,~$M_{k,C}(x)$ is a principal submatrix of~$M_{k,D}(x)$ and hence positive semidefiniteness of~$M_{k,D}(x)$ implies positive semidefiniteness of~$M_{k,C}(x)$. For computing~$B_4(n,d,w)$, it suffices to require that the matrices~$M_{3,D}(x)$ for each~$D \in \mathcal{C}_{k-1}$ with~$|D| \leq 1$ and the matrices~$M_{4,D}(x)$ for each~$D \in \mathcal{C}_{k}$ with~$|D|\geq 2$ even are positive semidefinite.

If~$|D|=k$, then~$M_{k,D}(x)$ is a matrix of order~$1 \times 1$, so positive semidefiniteness of~$M_{k,D}(x)$ is equivalent to~$x(D) \geq 0$.
We can assume in~$(\ref{A4ndw})$ and~$(\ref{Bndw})$ that~$x \, : \, \mathcal{C}_k \to \mathbb{R}_{\geq 0}$, since if~$|D|\leq k$ then~$x(D)$ occurs on the diagonal of~$M_{k,D}(x)$ and if~$|D| < 2$ and~$k \geq 3$ then~$x(D)$ occurs on the diagonal of~$M_{k-1,D}(x)$.

\subsection{Exploiting the symmetry of the problem}
 Fix~$k \in \N$ with~$k \geq 2$. The group~$G:=S_n$ acts naturally on~$\mathcal{C}_k$ by simultaneously permuting the indices~$1,\ldots,n$ of each code word in~$C \in \mathcal{C}_k$ (since the weight of each codeword is invariant under this action), and this action maintains minimum distances and cardinalities of codes~$C \in \mathcal{C}_k$. We can assume that the optimum~$x$ in~$(\ref{A4ndw})$ (or~$(\ref{Bndw})$) is $G$-invariant, i.e., $x \circ g = x$ for all~$g \in G$. To see this, let~$x$ be an optimum solution for~$(\ref{A4ndw})$ (or~$(\ref{Bndw})$). For each~$g \in G$, the function~$x \circ g$ is again an optimum solution, since the objective value of~$x \circ g$ equals the objective value of~$x$ and~$x \circ g$ still satisfies all constraints in~$(\ref{A4ndw})$ (or~$(\ref{Bndw})$). Since the feasible region is convex, the optimum~$x$ can be replaced by the average of~$x \circ g$ over all~$g \in G$. This yields a $G$-invariant optimum solution.  
 
 Let~$\Omega_k$ be the set of~$G$-orbits on~$\mathcal{C}_k$. Then~$|\Omega_k|$ is bounded by a polynomial in~$n$. Since there exists a~$G$-invariant optimum solution, we can replace, for each~$\omega \in \Omega_k$ and~$C \in \omega$, each variable~$x(C)$ by a variable~$y(\omega)$. Hence, the matrices~$M_{j,D}(x)$ become matrices~$M_{j,D}(y)$ and we have considerably reduced the number of variables in~$(\ref{A4ndw})$ and~$(\ref{Bndw})$.

It is only required that we check positive semidefiniteness of~$M_{j,D}(y)$ for one code~$D$ in each~$G$-orbit of~$\mathcal{C}_k$, as for each~$g \in G$, the matrix~$M_{j,g(D)}(y)$ can be obtained by simultaneously permuting rows and columns of~$M_{j,D}(y)$. We will describe how to reduce these matrices in size.

\begin{table}[H] \small 
\begin{center}
    \begin{tabular}{| r | r | r|| r |>{\bfseries}r |r | r| r|}
    \hline
    $n$ & $d$ & $w$ &  \multicolumn{1}{b{17mm}|}{best lower bound known}  & \multicolumn{1}{b{21mm}|}{\textbf{new upper bound}} & \multicolumn{1}{b{18mm}|}{best upper bound previously known} & \multicolumn{1}{b{18.1mm}|}{$\floor{A_3(n,d,w)}$}& \multicolumn{1}{b{13mm}|}{Delsarte bound}\\\hline 
    17 & 6 & 7  & 166 & 206*  & 207 & 228 & 249 \\ 
    18 & 6 & 7  & 243 & 312*  & 318 & 353 & 408 \\ 
    19 & 6 & 7  & 338 & 463*  & 503 & 526 & 553\\ 
    19 & 6 & 8  & 408 & 693  & 718 & 718 & 751\\ 
    20 & 6 & 8  & 588 & 1084  &1106 &1136 & 1199\\  
    21 & 6 & 8  & 775 & 1665 &   1695 & 1772 & 1938\\
    21 & 6 & 9  & 1186& 2328 & 2359   &2359 & 2364\\     \hline 
    25 & 8 & 8  & 759 & 850 &  856  & 926  & 948 \\     
    21 & 8 & 9  & 280 & 294 &   302 & 314 & 358\\ 
    22 & 8 & 9  & 280 & 440 &  473 & 473 & 597 \\ 
    23 & 8 & 9  & 400 & 662 &  703  & 707  & 830\\
    24 & 8 & 9  & 640 & 968 &  1041 & 1041  & 1160\\
    25 & 8 & 9  & 829 & 1366 &  1486 & 1486 & 1626 \\
    26 & 8 & 9  & 887 & 1901 &  2104 &  2108 & 2282 \\
    27 & 8 & 9  & 1023 & 2616 & 2882   & 2918 & 3203 \\  
    22 & 8 & 10 & {\color{donkerrood}616} & {\color{donkerrood}616} &630   & 634 & 758 \\
    22 & 8 & 11 & {\color{donkerrood}672} & {\color{donkerrood}672} &  680  & 680 & 805\\    \hline
  %     24 & 10 & 9  &  56& 116 &   118 &118 & 119 al een bekende bovengrens, zie artikel etzion en co over linear programming bounds
    27 & 10 & 9  &  118 & 291 &  293  &299 & 299  \\ 
    22 & 10 & 10  & 46 & 71 &   72  &  72 & 82 \\
    23 & 10 & 10  & 54 & 116 &  117 &117 & 117 \\   
    26 & 10 & 10  & 130 &  397&  406& 406  & 412 \\    
    27 & 10 & 10  & 162 & 555 & 571  & 571& 579  \\    
    22 & 10 & 11  & 46 & 79 &  80  &80 &88  \\\hline
    \end{tabular}
\end{center}
  \caption{\small An overview of the new upper bounds for constant weight codes. The best previously known upper bounds are taken from Brouwer's table~$\cite{brouwertable}$. The unmarked new upper bounds are instances of~$B_4(n,d,w)$, and the new upper bounds marked with an asterisk are instances of~$A_4(n,d,w)$ (rounded down). For comparison: $B_4(17,6,7)=213$, $B_4(18,6,7)=323$ and~$B_4(19,6,7)=486$. The best previously known upper bounds on~$A(22,8,10)$ and~$A(22,8,11)$ were found in~$\cite{KT}$ and~$\cite{schrijver}$, respectively. It can be concluded that~$A(22,8,10)=616$ and~$A(22,8,11)=672$ (marked in red).\label{28}}
\end{table}

 For~$D \in \mathcal{C}_k$, let~$G_D$ be the subgroup of~$G$ consisting of all~$g \in G$ that leave~$D$ invariant. Then the action of~$G$ on~$\mathcal{C}_k$ induces an action of $G_D$ on~$\mathcal{C}_j(D)$.  The simultaneous action of~$G_D$ on the rows and columns of~$M_{j,D}(y)$ leaves~$M_{j,D}(y)$ invariant. Therefore, there exists a block-diagonalization~$M_{j,D}(y) \mapsto U^T  M_{j,D}(y) U$ of~$M_{j,D}(y)$, for a matrix~$U$ not depending on~$y$. Then~$M_{j,D}(y)$ is positive semidefinite if and only if each of the blocks is positive semidefinite. There are several equal (or equivalent) blocks and after removing duplicate (or equivalent) blocks we obtain a matrix of order bounded polynomially in~$n$ where the entries in each block are linear functions in the variables~$y(\omega)$ (with coefficients bounded polynomially in~$n$). Hence, we have reduced the size of the matrices involved in our semidefinite program.

Note that, after reductions, the number of variables involved in the semidefinite programs for computing~$A_k(n,d,w)$ and~$B_k(n,d,w)$ are the same. However, the program for computing~$B_k(n,d,w)$ contains fewer blocks than the program for computing~$A_k(n,d,w)$, and the blocks are smaller and contain fewer variables per matrix entry. This is important, as the semidefinite programs for computing~$A_4(n,d,w)$ turn out to be very large in practice (although they are of polynomial size).

For particular weights~$w$ (in the case of constant weight codes), the group of distance preserving permutations of~$\mathcal{C}_k$ can be larger than~$S_n$. If~$w=n/2$ there is a further action of~$S_2$ on~$\mathcal{C}_k$ by adding the all-ones word to each word in each code in~$\mathcal{C}_k$. Since the corresponding reduction of the semidefinite program can only be used for specific weights~$w$, we do not consider the reduction in this paper, although it was used for reducing the number of variables in computing~$B_4(22,8,11)$.

The reductions of the optimization problem will be described in detail in Section~$\ref{red}$. Table~$\ref{28}$ contains the new upper bounds for~$n\leq 28$, which are the values of~$n$ usually considered. Since some tables on Brouwer's website~$\cite{brouwertable}$ also consider~$n$ in the range~$29 \leq n \leq 32$, we give some new bounds for these cases (many of which are computed with the smaller program~$A_3(n,d,w)$) in Table~$\ref{32}$. The paper is concluded by two appendices, one giving pseudocode (which outlines the structure of the program for generating the semidefinite program) and one specifiying a subroutine for computing the polynomials involved. All improvements have been found using multiple precision versions of SDPA~\cite{sdpa}, where the largest program (for computing~$B_4(22,8,10)$) took approximately three weeks to compute on a modern desktop pc. 

\subsection{Comparison with earlier bounds}

In this paper we will consider~$B_4(n,d,w)$ and~$A_4(n,d,w)$. The reduction, based on representation theory, is an adaptation to constant weight codes of the method in~$\cite{onsartikel}$. The method of~$\cite{onsartikel}$ was also used (for~$k=3$) in the context of mixed binary/ternary codes~$\cite{mixed}$. It can be proved that~$A_2(n,d,w)$ is equal to the Delsarte bound~$\cite{delsarte}$. The bound~$A_4(n,d,w)$ is an adaptation of the bound~$A_4(n,d)$ for non-constant weight codes considered in~$\cite{semidef}$. The semidefinite programming bound for constant weight codes introduced by Schrijver in~$\cite{schrijver}$ is a slight sharpening of~$A_3(n,d,w)$ (in almost all cases it is equal to~$A_3(n,d,w)$).  The bound~$B_4(n,d,w)$, which is based on quadruples of code words, is a bound `in between'~$A_3(n,d,w)$ and~$A_4(n,d,w)$: it is the bound~$A_3(n,d,w)$ with constraints for matrices~$M_{4,D}(x)$ with~$|D|=2$ (based on quadruples of code words) added. Or it can be seen as a bound obtained from~$A_4(n,d,w)$ by removing the (large) matrix~$M_{4,\emptyset}(x)$ and replacing it by~$M_{3,\emptyset}(x)=M_{2,\emptyset}(x)$. 

Recently, Kim and Toan~$\cite{KT}$ added linear inequalities to the Schrijver bound~$\cite{schrijver}$.  An advantage of their method is that the semidefinite programs are small and can be solved fast.  The bound~$B_4(n,d,w)$ is often sharper than their bound, but it takes much more time to compute. Only the cases of~$n,d,w$ with~$n\leq 28$ for which finding~$B_4(n,d,w)$  did not require excessive computing time or memory are therefore considered in the present work.

\section{Preliminaries on representation theory}
In this section we give the definitions and notation from representation theory (where we are mostly concerned with the symmetric group) used throughout the paper, similarly to the notation used in~$\cite{onsartikel}$. Proofs are omitted, but for more information, the reader can consult Sagan~$\cite{sagan}$. 

A \emph{group action} of a group~$G$ on a set~$X$ is a group homomorphism~$\phi: G \to S_X$, where~$S_X$ is the group of bijections of~$X$ to itself. If~$G$ acts on~$X$, we write~$g \circ x:= \phi(g)(x)$ for all~$g \in G$ and~$x \in X$ and we write~$X^G$ for the set of elements of~$X$ invariant under the action of~$G$. If~$X$ is a linear space, the elements of~$S_X$ are assumed to be linear functions. The action of~$G$ on a set~$X$   induces an action of~$G$ on the linear space~$\mathbb{C}^X$, by~$(g \circ f)(x):= f(g^{-1} \circ x)$, for~$g\in G$,~$f \in \mathbb{C}^X$ and~$x \in X$.    

If~$m \in \N$ and~$G$ is a finite group acting on~$V=\C^m$, then~$V$ is a~$G$\emph{-module}. If~$V$ and~$W$ are~$G$-modules, then a $G$\emph{-homomorphism} $\psi: V \to W$ is a linear map such that~$g \circ \psi(v)=\psi(g \circ v)$ for all~$g \in G$,~$v \in V$. Moreover, a module~$V$ is called \emph{irreducible} if the only~$G$-invariant submodules of~$V$ are~$\{0\}$ and~$V$ itself.

Suppose that~$G$ is a finite group acting \emph{unitarily} on~$V=\C^m$. This means that for each~$g\in G$ there is a unitary matrix~$U_g \in \C^{m \times m}$ such that~$g \circ x = U_gx$ for all~$x \in \C^m$. Consider the inner product~$\langle x, y \rangle :=x^*y$ for~$x,y \in \C^m$, where~$x^*$ denotes the conjugate transpose of~$x \in \C^m$. Then~$V$ can be decomposed as a direct sum of~$G$\emph{-isotypical components}~$V_1,\ldots,V_k$. This means that~$V_i$ and~$V_j$ are orthogonal for distinct~$i$ and~$j$ (with respect to the mentioned inner product), and each~$V_i$ is a direct sum~$V_{i,1} \oplus \ldots \oplus V_{i,m_i}$ of irreducible and mutually isomorphic~$G$-modules, such that~$V_{i,j}$ and~$V_{i,j'}$ are orthogonal for distinct~$j, j'$ and such that~$V_{i,j}$ and~$V_{i',j'}$ are isomorphic if and only if~$i=i'$.  

For each~$i \leq k$ and~$j \leq m_i$ we choose a nonzero vector~$u_{i,j} \in V_{i,j}$ with the property that for each~$i$ and all~$j,j'\leq m_i$ there exists a~$G$-isomorphism~$V_{i,j} \to V_{i,j'}$ mapping~$u_{i,j}$ to~$u_{i,j'}$. For each~$i \leq k$, we define~$U_i$ to be the matrix~$[u_{i,1},\ldots ,u_{i,m_i}]$ with columns~$u_{i,j}$ ($j=1,\ldots,m_i$). Any set of matrices~$\{U_1,\ldots,U_k\}$ obtained in this way is called a \emph{representative} set for the action of~$G$ on~$\C^m$. Then the map
\begin{align} \label{cPhi}
    \Phi : (\C^{m \times m})^G \to \bigoplus_{i=1}^k \C^{m_i \times m_i} \,\, \text{ with } \,\, A \mapsto \bigoplus_{i=1}^k U_i^* A U_i
\end{align}
is bijective. So~$ \dim((\C^{m \times m})^G) =\sum_{i=1}^k m_i^2$, which can be considerably smaller than~$m$. Another crucial property for our purposes is that any~$A \in (\C^{m \times m})^G $ is positive semidefinite (i.e., self-adjoint with all eigenvalues nonnegative) if and only if the image~$\Phi(A)$ is positive semidefinite. 

In this paper,~$G$ is acting \emph{real-orthogonally} on a vector space~$V=\R^m$. This means that for each~$g \in G$ there is a real orthogonal matrix~$U_g \in \R^{m \times m}$ with~$g \circ x = U_g x$ for all~$x \in \C^m$. Moreover, it turns out that all representative sets we define consist of real matrices. Then
\begin{align} \label{PhiR} 
\Phi(A) = \bigoplus_{i=1}^k U_i^T A U_i \text{ for }A \in (\R^{m \times m})^G,
\end{align}
and~$\Phi\left((\R^{m \times m})^G\right) = \oplus_{i=1}^k \R^{m_i \times m_i}$.
Also,~$A \in \R^{ m \times m}$ is positive semidefinite if and only if each of the matrices~$U_i^T A U_i$ is positive semidefinite ($i=1,\ldots,k$). This is very useful for checking whether~$A$ is positive semidefinite. Since~$V_{i,j}$ is the linear space spanned by~$G \circ u_{i,j}$ (for each~$i,j$), we have~$\R^m = \oplus_{i=1}^k\oplus_{j=1}^{m_i} \R G \circ u_{i,j}$,
where~$\R G$ denotes the group algebra of~$G$. It will be convenient to consider the columns of~$U_i$ as elements of the dual space~$(\R^m)^*$ via the inner product mentioned above.

\subsection{A representative set for the action of~\texorpdfstring{$S_n$}{Sn} on~\texorpdfstring{$V^{\otimes n}$}{Vtensorn}}\label{repr}
Fix~$n\in \N$. We will consider the natural action of~$S_n$ on~$V^{\otimes n}$, where~$V$ is a finite-dimensional vector space, by permuting the indices. We describe a representative set for the action of~$S_n$ on~$V^{\otimes n}$, that will be used repeatedly in the reductions throughout this paper. 

A \emph{partition}~$\lambda $ of~$n$ is a sequence~$(\lambda_1,\ldots, \lambda_h)$ of natural numbers with~$\lambda_1 \geq \ldots \geq \lambda_h >0$ and~$\lambda_1 + \ldots + \lambda_h = n$. The number~$h$ is called the \emph{height} of~$\lambda$. We write~$\lambda \vdash n$ if~$\lambda$ is a partition of~$n$. The~\emph{Young shape} (or \emph{Ferrers diagram})~$Y(\lambda)$ of~$\lambda$ is the set
\begin{align}
    Y(\lambda) := \{(i,j) \in \N^2 \, | \, 1 \leq j \leq h, \, 1 \leq i \leq \lambda_j\}.  
\end{align}
Fixing an index~$j_0 \leq h$, the set of elements~$(i,j_0)$ (for~$1 \leq i \leq \lambda_j$) in~$Y(\lambda)$ is called the~$j_0$\emph{-th row} of~$Y(\lambda)$. Similarly, fixing an element~$i_0 \leq \lambda_1$, the set of elements~$(i_0,j)$ (where~$j$ varies) in~$Y(\lambda)$ is called the~$i_0$\emph{-th column} of~$Y(\lambda)$.  Then the \emph{row stabilizer}~$R_{\lambda}$ of~$\lambda$ is the group of permutations~$\pi$ of~$Y(\lambda)$ with~$\pi(Z)=Z$ for each row~$Z$ of~$Y(\lambda)$ and, similarly, the \emph{column stablizer}~$C_{\lambda}$ of~$\lambda$ is the group of permutations~$\pi$ of~$Y(\lambda)$ with~$\pi(Z)=Z$ for each column~$Z$ of~$Y(\lambda)$. 

A \emph{Young tableau} with shape~$\lambda$ (also called a~$\lambda$\emph{-tableau}) is a function~$\tau : Y(\lambda) \to \N$. A Young tableau with shape~$\lambda$ is \emph{semistandard} if the entries are nondecreasing in each row and strictly increasing in each column. Let~$T_{\lambda,m}$ be the collection of semistandard $\lambda$-tableaux with entries in~$[m]$. Then~$T_{\lambda,m} \neq \emptyset$  if and only if~$m$ is at least the height of~$\lambda$. We write~$\tau \sim \tau'$ for~$\lambda$-tableaux~$\tau,\tau´$ if~$\tau'=\tau r$ for some~$r \in R_{\lambda}$.

Let~$B=(B(1),\ldots,B(m))$ be an ordered basis of~$V^*$. For any~$\tau \in T_{\lambda,m}$, define 
\begin{align} \label{utau}
    u_{\tau,B}:= \sum_{\tau'\sim \tau } \sum_{c \in C_{\lambda}} \text{sgn}(c) \bigotimes_{y \in Y(\lambda)} B\left(\tau '(c(y))\right).
\end{align}
Here the Young shape~$Y(\lambda)$ is ordered by concatenating its rows. Then (cf.~$\cite{sagan}$ and~$\cite{onsartikel}$) the set 
\begin{align} \label{reprsetdef}
\left\{\, [u_{\tau,B} \,\, | \,\, \tau \in  T_{\lambda,m}] \,\, | \,\, \lambda \vdash n  \right\},
\end{align}
consisting of matrices, is a representative set for the natural action of~$S_n$ on~$V^{\otimes n}$. 

\section{Reduction of the optimization problem}\label{red}
In this section we give the reduction of optimization problem~$(\ref{Bndw})$ for computing~$B_4(n,d,w)$, using the representation theory from the previous section. Also, we give a reduction for computing~$A_4(n,d,w)$. First we consider block diagonalizing~$M_{4,D}(y)$ for~$D \in \mathcal{C}_4$ with~$|D|=1$ or~$|D|=2$, for computing both~$B_4(n,d,w)$ and~$A_4(n,d,w)$.\footnote{Note that if~$|D|=1$ then~$M_{4,D}(y)=M_{3,D}(y)$. For computing~$A_4(n,d,w)$, we only need to consider the case~$|D|=2$ as~$M_{4,D}(y)$ for~$|D|=1$ is a principal submatrix of~$M_{4,\emptyset}(y)$.}  Subsequently we consider the cases~$M_{3,\emptyset}(y) =M_{2,\emptyset}(y)$ or~$M_{4,\emptyset}(y)$, for computing~$B_4(n,d,w)$ or~$A_4(n,d,w)$, respectively. Note that for the cases~$|D|=3$ and~$|D|=4$ the matrix~$M_{4,D}(y) = (y(D))$  is its own block diagonalization, so then~$M_{4,D}(y)$ is positive semidefinite if and only if~$y(D) \geq 0$. 

\subsection{The cases~\texorpdfstring{$|D|=1$}{|D|=1} and~\texorpdfstring{$|D|=2$}{|D|=2}} \label{D2}
In this section we consider one code~$D \in \mathcal{C}_4$  with~$|D|=1$ or~$|D|=2$. We can assume that~$D=\{v_1,v_2\}$ with
\begin{align}
    v_1=&\overbrace{1\ldots1\,1\ldots1}^{w}\overbrace{0\ldots0\,\,0\ldots0}^{n-w} 
    \\v_2= &\underbrace{0\ldots0}_{t}\underbrace{1\ldots1\,1\ldots1}_{w}\underbrace{\,0\ldots0}_{n-t-w},    \notag 
\end{align}
where~$t \in \Z_{\geq 0}$ with~$t=0$ or~$d/2 \leq t \leq w$. For the remainder of this section, fix~$t \in \{0\} \,\cup\, \{t \,\,|\,\, d/2 \leq t \leq w \}$ (recall that~$w \leq n/2$, so~$t \leq n-w$). If~$t=0$, then~$|D|=1$ and if~$d/2 \leq t \leq w$, then~$|D|=2$. The rows and the columns of~$M_{4,D}(y)$ are parametrized by codes~$C\supseteq D$ of size at most~$3$ (if~$|D|=2$) or size at most~$2$ (if~$|D|=1$). 

Let~$H$ be the group of distance preserving permutations of~$\mathcal{C}_4$ that fix~$v_1$ and~$v_2$. So
\begin{align}
    H \cong S_t \times S_{w-t} \times S_t \times S_{n-t-w}.
\end{align}
We first describe a representative set for the action of~$H$ on~$\mathbb{R}^{\mathbb{F}_2^n}$ and then restrict to words of weight~$w$ and distance at least~$d$ to both words in~$D$. Let~$e_j$ denote the~$j$-th standard basis vector of~$\mathbb{R}^{\mathbb{F}_2}$, for~$j=0,1$. Define~$B=(B(1),B(2)):=(e_0, e_1)$, where we consider~$B(1)$ and~$B(2)$ as elements of the dual space~$(\mathbb{R}^{\mathbb{F}_2})^*$ via the standard inner product. 

Fix~$\bm{n}=(n_1,n_2,n_3,n_4):=(t,w-t,t,n-t-w)$ and let~$\bm{\lambda} \vdash \bm{n}$ mean that~$\bm{\lambda}=(\lambda_1,\ldots,\lambda_4)$ with~$\lambda_i \vdash n_i$ for~$i=1,\ldots,4$ (i.e.,\ each~$\lambda_i$ is equal to~$(\lambda_{i,1} ,\ldots , \lambda_{i,h})$ for some~$h$). For~$\bm{\lambda} \vdash \bm{n}$, define
\begin{align}
    W_{\bm{\lambda}} := T_{\lambda_1,2} \times T_{\lambda_2,2} \times T_{\lambda_3,2} \times T_{\lambda_4,2},
\end{align}
and for~$\bm{\tau}=(\tau_1,\ldots,\tau_4) \in W_{\bm{\lambda}}$, define (cf.~(\ref{utau}))
\begin{align} 
\bm{u_{\bm{\tau}}} := \bigotimes_{i=1}^4 u_{\tau_i, B}.
\end{align}
Now a representative set for the action of~$H$ on~$\mathbb{R}^{\mathbb{F}_2^n}$, using the natural isomorphism $ \mathbb{R}^{\mathbb{F}_2^n} \cong (\mathbb{R}^{\mathbb{F}_2})^{\otimes n}$, is (cf.\ Section~$\ref{repr}$)
\begin{align}
\{[\bm{u_{\bm{\tau}}}  \,\, | \,\, \bm{\tau} \in      W_{\bm{\lambda}}] \,\, | \,\, \bm{\lambda} \vdash \bm{n}   \}. 
\end{align}
We restrict to words of weight~$w$ and distance contained in~$\{0,d,d+1,\ldots,n\}$ to both words in~$D$. For~$d_1,d_2\in \{0,1,\ldots,n\}$, let~$L_{w,d_1,d_2}$ denote the linear subspace of~$\mathbb{R}^{\mathbb{F}_2^n}$ spanned by the unit vectors~$\epsilon_v$, with~$v$ a word of weight~$w$ and distances~$d_1$ and~$d_2$ to~$v_1$ and~$v_2$, respectively. Then~$L_{w,d_1,d_2}$ is~$H$-invariant. Moreover, for any~$\bm{\lambda} \vdash \bm{n}$ and~$\bm{\tau}=(\tau_1,\ldots,\tau_4) \in W_{\bm{\lambda}}$, the irreducible representation~$\R H \cdot \bm{u_{\bm{\tau}}}$ is contained in~$L_{w,d_1,d_2}$, with
\begin{align}
  w&=  |\tau_1^{-1}(2)|+|\tau_2^{-1}(2)|+|\tau_3^{-1}(2)|+|\tau_4^{-1}(2)|,\\ d_1 &=|\tau_1^{-1}(1)|+|\tau_2^{-1}(1)|+|\tau_3^{-1}(2)|+|\tau_4^{-1}(2)|, \notag  \\ 
    d_2&=  |\tau_1^{-1}(2)|+|\tau_2^{-1}(1)|+|\tau_3^{-1}(1)|+|\tau_4^{-1}(2)|.\notag 
\end{align}
Let 
\begin{align} \label{wlambda}
    W_{\bm{\lambda}}' := \{\bm{\tau} \in W_{\bm{\lambda}} \, | \,\,\, & |\tau_1^{-1}(2)|+|\tau_2^{-1}(2)|+|\tau_3^{-1}(2)|+|\tau_4^{-1}(2)|=w, \\ & |\tau_1^{-1}(1)|+|\tau_2^{-1}(1)|+|\tau_3^{-1}(2)|+|\tau_4^{-1}(2)| \in \{0,d,d+1,\ldots,n\}, \notag  \\ 
    &  |\tau_1^{-1}(2)|+|\tau_2^{-1}(1)|+|\tau_3^{-1}(1)|+|\tau_4^{-1}(2)| \in \{0,d,d+1,\ldots,n\} \}. \notag
\end{align}
Furthermore, let~$Z$ be the~$\mathcal{C}_4(D) \times \mathbb{F}_2^n$ matrix with~$0,1$ entries satisfying
\begin{align}
    Z_{C,\alpha} = 1 \,\, \text{ if and only if } C =   \{v_1,v_2,\alpha\},
\end{align}
for~$C \in \mathcal{C}_4(D)$ and~$\alpha \in \mathbb{F}_2^n$. The map~$x \mapsto Zx$ is a surjective~$H$-homomorphism~$\mathbb{R}^{\mathbb{F}_2^n} \to \mathbb{R}^{\mathcal{C}_4(D)}$. 
Then 
\begin{align} \label{set}
\{ZU_{\bm{\lambda}}\,\, | \,\, \bm{\lambda} \vdash \bm{n}   \}, \,\,\,\text{ where }\,\, U_{\bm{\lambda}} := [\bm{u_{\bm{\tau}}}  \,\, | \,\, \bm{\tau} \in      W_{\bm{\lambda}}']  \,\, \text{ (for $\bm{\lambda} \vdash \bm{n}$)}, 
\end{align}
  is a representative set for the action of~$H$ on~$\R^{\mathcal{C}_4^d(D)} \subseteq \R^{\mathcal{C}_4(D)} $, where~$\mathcal{C}_4^d(D)$ denotes the collection of codes~$C \in \mathcal{C}_4(D)$ with~$d_{\text{min}}(C) \geq d$.

\subsubsection{Computations for~\texorpdfstring{$|D|=1$ or~$|D|=2$}{D=1 or D=2}} \label{d12comp}
Fix~$D = \{v_1,v_2\} \in \mathcal{C}_4$. Let~$\Omega_4(D)$ denote the set of all~$S_n$-orbits of codes in~$\mathcal{C}_4$ containing~$D=\{v_1,v_2\}$. For each~$\omega \in \Omega_4(D)$, we define the~$\mathcal{C}_4(D) \times \mathcal{C}_4(D)$-matrix~$N_{\omega}$ with entries in~$\{0,1\}$ by
\begin{align}
    (N_{\omega})_{\{v_1,v_2,\alpha\},\{v_1,v_2,\beta\}} := \begin{cases} 1 &\mbox{if } \{v_1,v_2,\alpha,\beta\} \in \omega,  \\ 
0 & \mbox{else.} \end{cases} 
\end{align}
Then we obtain with~$(\ref{PhiR})$ and~$(\ref{set})$ that, for each~$y : \Omega_4(D) \to \mathbb{R}$,
\begin{align} \label{blocks1t}
    \Phi(M_{4,D}(y)) = \Phi \left(\sum_{\omega \in \Omega_4(D)}y(\omega) N_{\omega} \right ) = \bigoplus_{\bm{\lambda} \vdash \bm{n}} \sum_{\omega \in \Omega_4(D)} y(\omega) (ZU_{\bm{\lambda}})^T N_{\omega} ZU_{\bm{\lambda}}. 
\end{align}
Define~$N_{\omega}':=Z^TN_{\omega}Z$, for any~$\omega \in \Omega_4(D)$. So~$N_{\omega}'$ is an~$\mathbb{F}_2^n \times \mathbb{F}_2^n$-matrix with
\begin{align}
    (N_{\omega}')_{\alpha,\beta} := \begin{cases} 1 &\mbox{if } \{v_1,v_2,\alpha,\beta\} \in \omega,  \\ 
0 & \mbox{else.} \end{cases} 
\end{align}
Then~$(\ref{blocks1t})$ becomes
\begin{align} \label{blocks1}
    \Phi(M_{4,D}(y))  = \bigoplus_{\bm{\lambda} \vdash \bm{n}} \sum_{\omega \in \Omega_4(D)} y(\omega) U_{\bm{\lambda}}^T N_{\omega}' U_{\bm{\lambda}}. 
\end{align}
The number of~$\bm{\lambda} \vdash \bm{n}$, and the numbers~$|W_{\bm{\lambda}}'|$ and~$|\Omega_4(D)|$  are all polynomially bounded in~$n$. Hence the number of blocks in~$(\ref{blocks1})$, as well as the size of each block and the number of variables occuring in all blocks are bounded by a polynomial in~$n$. In the remainder of this section we show how to compute the coefficients in the blocks~$U_{\bm{\lambda}}^T N_{\omega}' U_{\bm{\lambda}}$ in polynomial time, that is, we will show how to compute each entry~$\bm{u_{\bm{\tau}}}^T N_{\omega}' \bm{u_{\bm{\sigma}}}$ in polynomial time, for~$\bm{\tau}, \bm{\sigma} \in W_{\bm{\lambda}}'$. 

For~$P=(i_1,i_2,i_3,i_4) \in \mathbb{F}_2^{4}$, let~$a_P := e_{i_1} \otimes e_{i_2} \otimes e_{i_3} \otimes e_{i_4}$,
where~$e_{j}$ denotes the standard unit basis vector in~$\R^{\mathbb{F}_2}$ corresponding to~$j \in \mathbb{F}_2$. Then the set~$E_4:=\{a_P \, | \, P \in \mathbb{F}_2^{4} \}$ is a basis of~$(\R^{\mathbb{F}_2})^{\otimes 4}$ and we  define~$E_4^*$ to be the dual basis. Similarly, for~$P'=(i_1,i_2) \in \mathbb{F}_2^{2}$, let~$a_{P'} := e_{i_1} \otimes e_{i_2}$. Then the set~$E_2:=\{a_{P'} \, | \, P' \in \mathbb{F}_2^{2} \}$ is a basis of~$(\R^{\mathbb{F}_2})^{\otimes 2}$ and we  define~$E_2^*$ to be the dual basis. Furthermore, for any~$P=(i_1,i_2,i_3,i_4) \in \mathbb{F}_2^{4}$ we define~$h(P):=(i_3,i_4) \in \mathbb{F}_2^{2} $. 

For any monomial~$p=a_{P_1}^*\cdot\ldots\cdot a_{P_n}^*$ of degree~$n$ on~$(\R^{\mathbb{F}_2})^{\otimes 4}$ (with~$P_1,\ldots,P_n \in \mathbb{F}_2^{4}$), we define, for~$(i_1,i_2) \in \mathbb{F}_2^2$: 
\begin{align}
 \xi_{i_1,i_2}(p) = \xi_{i_1,i_2}(a_{P_1}^*\cdot\ldots\cdot a_{P_n}^*)  := \prod_{\substack{i \in \{1,\ldots,n\}: \\ P_i \in \{i_1\} \times \{i_2\} \times \mathbb{F}_2^2 }} a_{h(P_i)}^*,
\end{align}
which is a monomial on~$(\R^{\mathbb{F}_2})^{\otimes 2}$.

Let~$Q$ be the set of monomials~$p$ of degree~$n$ on~$(\R^{\mathbb{F}_2})^{\otimes 4}$, that satisfy~$\deg(\xi_{1,0}(p))=\deg(\xi_{0,1}(p))=t$, $\deg(\xi_{1,1}(p))=w-t$ and~$\deg(\xi_{0,0}(p))=n-w-t$. For any~$(\alpha,\beta) \in (\F_2^n)^2$, define the following element of~$Q$:
\begin{align}\label{agebeuren}
 \psi(\alpha,\beta) :=  \prod_{i=1}^t a_{(1,0,\alpha_i,\beta_i)}^* \cdot \prod_{i=t+1}^{w} a_{(1,1,\alpha_i,\beta_i)}^* \cdot \prod_{i=w+1}^{w+t} a_{(0,1,\alpha_i,\beta_i)}^* \cdot  \prod_{i=w+t+1}^n a_{(0,0,\alpha_i,\beta_i)}^*.
\end{align}
Then~$\psi(\alpha,\beta)=\psi(\alpha',\beta')$ if and only if~$(\alpha,\beta)$ and~$(\alpha',\beta')$ belong to the same~$H$-orbit on~$(\mathbb{F}_2^n)^2$.  Hence,~$(\ref{agebeuren})$ gives a bijection between the set~$Q$ and the set of~$H$-orbits on~$(\mathbb{F}_2^n)^2$. Write~$\mathcal{C}_4'$ for all codes~$C \subseteq \mathbb{F}_2^n$ of size~$\leq 4$ (so not necessarily of constant weight~$w$). Then the function
\begin{align}
    \left( \mathbb{F}_2^t \times \mathbb{F}_2^{w-t} \times \mathbb{F}_2^t \times  \mathbb{F}_2^{n-w-t} \right)^2 &\to \mathcal{C}_4',\\
((\alpha^1,\alpha^2,\alpha^3,\alpha^4),(\beta^1,\beta^2,\beta^3,\beta^4)) &\mapsto \{v_1,v_2,\alpha^1\alpha^2\alpha^3\alpha^4, \beta^1\beta^2\beta^3\beta^4\},    \notag 
\end{align}
induces a surjective function~$  r\, : \, Q \to \Omega_4'(D)$, where~$\Omega_4'(D)\supseteq \Omega_4(D)$ denotes the set of all~$H$-orbits of codes in~$\mathcal{C}_4'$ that contain~$D=\{v_1,v_2\}$. 

For any $\mu \in Q$, define
\begin{align} 
A_{\mu}:= \Big\{(a_i)_{i=1}^n \in E_2^t \times E_2^{w-t} \times E_2^{t}\times E_2^{n-w-t}\,\, \Big\rvert &\,\,\,\prod_{i=1}^t a_i^* = \xi_{1,0}(\mu), \, \prod_{i=t+1}^w a_i^* = \xi_{1,1}(\mu), 
\\ &\prod_{i=w+1}^{w+t} a_i^* = \xi_{0,1}(\mu), \,\prod_{i=w+t+1}^n a_i^* = \xi_{0,0}(\mu) \Big\}, \notag
\end{align} 
and 
\begin{align}
K_{\mu} := \sum_{(a_1,\ldots,a_n) \in A_{\mu} }
 \left( \bigotimes_{i=1}^t a_i\right) \otimes \left( \bigotimes_{i=t+1}^w a_i\right) \otimes  \left( \bigotimes_{i=w+1}^{w+t} a_i\right) \otimes  \left( \bigotimes_{i=w+t+1}^n a_i\right). 
\end{align}
(Here every~$a_i \in (\mathbb{R}^{\F_2})^{\otimes 2}$ is regarded as an element of~$\mathbb{R}^{\F_2 \times \F_2}$ via the natural isomorphism, so that~$K_{\mu}$ is an~$\mathbb{F}_2^n \times \mathbb{F}_2^n$-matrix.)
Then Lemma~$1$ of~$\cite{onsartikel}$ implies that, if~$\omega \in \Omega_4(D)$, then
\begin{align}
    N_{\omega}' = \sum_{\substack{\mu \in Q\\ r(\mu) = \omega  }} K_{\mu}.
\end{align}
Define the following elements of~$E_2^*$:
 \begin{align}
    \eta_1 &:= a_{(1,0)}^*,   \,\,\,\,\,\,     \eta_3 = a_{(0,1)}^*,  \\ 
    \eta_2 &:= a_{(1,1)}^*,   \,\,\,\,\,\,     \eta_4 = a_{(0,0)}^*. \notag 
 \end{align}
Furthermore, we define for~$\bm{\tau}$ and~$\bm{\sigma}$ in~$\bm{W}_{\bm{\lambda}}'$, the following polynomial on~$(\R^{\mathbb{F}_2})^{\otimes 4}$:
\begin{align} \label{firstp}
    p_{\bm{\tau},\bm{\sigma}} := \prod_{i=1}^4 \sum_{\substack{\tau_i'\sim \tau_i \\ \sigma_i' \sim \sigma_i  }} \sum_{c_i, c_i' \in C_{\lambda_i}} \text{sgn}(c_ic_i') \prod_{y \in Y(\lambda_i)} \eta_i \otimes \left( B(\tau_i'c_i(y)) \otimes B(\sigma_i'c_i'(y)) \right). 
\end{align}
 This polynomial can be computed (expressed as a linear combination of monomials in~$\eta_i \otimes ( B(j) \otimes B(h))$) in polynomial time, as proven in~$\cite{gijswijt, onsartikel}$. See the appendix for an algorithm due to Gijswijt~$\cite{gijswijt}$. Note that each monomial in~$ p_{\bm{\tau}, \bm{\sigma}}$ is contained in~$Q$. By Lemma~$2$ of~$\cite{onsartikel}$, we obtain
\begin{align}
 \sum_{\mu \in Q} (\bm{u_{\bm{\tau}}}^T K_{\mu}\bm{u_{\bm{\sigma}}})\mu = p_{\bm{\tau}, \bm{\sigma}},
\end{align}
which is a linear combination of~$\mu \in Q$. Hence one can compute the entry $\sum_{\omega \in \Omega_4(D)} y(\omega) \bm{u_{\bm{\tau}}}^T   N_{\omega}' \bm{u_{\bm{\sigma}}}$ by first expressing~$p_{\bm{\tau}, \bm{\sigma}}$ as a linear combination of~$\mu \in Q$ and subsequently replacing each $\mu \in Q$ in~$p_{\bm{\tau}, \bm{\sigma}}$ with the variable~$y(r(\mu))$ if~$r(\mu) \in \Omega_4(D)$ is an orbit of minimum distance~$\geq d$ and with zero otherwise.    

\subsection{The case~\texorpdfstring{$D=\emptyset$}{D=empty}} \label{Dempty}
Next, we consider how to block diagonalize~$M_{3,\emptyset}(y)=M_{2,\emptyset}(y)$ for computing~$B_4(n,d,w)$. Also we give a reduction of the matrix~$M_{4,\emptyset}(y)$ for computing~$A_4(n,d,w)$. So we will reduce the matrices~$M_{2s,\emptyset}(y)$ for~$s \in \{1,2\}$, where we consider~$s=1$ for computing~$B_4(n,d,w)$ and~$s=2$ for computing~$A_4(n,d,w)$. We start by giving a representative set for the natural action of~$S_n$ on~$(\R^{\F_2^s})^{\otimes n} \cong \R^{(\F_2^n)^s}$, using the results described in Section~\ref{repr}. 

Let~$J_s=(e_z \,\, | \,\, z \in \F_2^s)$ be an ordered~$2^s$-tuple containing the unit basis vectors of~$\R^{\F_2^s}$ as columns. Then we can view~$J_s$ as an ordered basis of~$(\R^{\F_2^s})^*$ via the standard inner product. So we have
\begin{align}\label{Cbasis}
    J_1 &= (J_1(1), J_1(2)) = (e_0,e_1),\\
    J_2 &= (J_2(1), J_2(2), J_2(3), J_2(4)) = (e_{0,0},e_{0,1},e_{1,0},e_{1,1}).\notag   
\end{align}
Then (cf.~$(\ref{reprsetdef})$)
\begin{align} \label{reprtot}
\{ [u_{\tau,J_s} \,\, | \,\, \tau \in  T_{\lambda,2^s}] \,\, | \,\, \lambda \vdash n  \}
\end{align}
is a representative set for the natural action of~$S_n$ on~$(\R^{\F_2^s})^{\otimes n}$.

For computing~$B_4(n,d,w)$, we first consider~$s=1$.  We restrict the representative set~$(\ref{reprtot})$ for the action of~$S_n$ on~$(\R^{\F_2})^{\otimes n} \cong \R^{\F_2^n}$ to~$\R ^F$ (recall:~$F \subseteq \F_2^n$ is the set of all words of constant weight~$w$). Let
\begin{align}
      R_{\lambda}^{(1)}:= \{\tau \in T_{\lambda, 2} \, | \,\,\,    |\tau^{-1}(2)|=w\}.
\end{align}
Then    
\begin{align} \label{set1}
\{U_{\lambda}^{(1)}:= [u_{\tau,J_1} \,\, | \,\, \tau \in  R_{\lambda}^{(1)}] \,\, | \,\, \lambda \vdash n  \}
\end{align}
is representative for the action of~$S_n$ on~$\R^{F}=\R^{\mathcal{C}_1 \setminus \{ \emptyset\}} \subseteq \R^{\F_2^n} $. 

For computing~$A_4(n,d,w)$, we consider~$s=2$. We proceed by restricting the representative set~$(\ref{reprtot})$ of the action of~$S_n$ on~$\R^{(\F_2^n)^2}$ to pairs of words in~$F^2$ with distance contained in~$\{0,d,d+1,\ldots,n\}$. Given~$w_1,w_2,d_1\in \{0,1,\ldots,n\}$, let~$L_{w_1,w_2,d_1}^{(2)}$ denote the linear subspace of~$\mathbb{R}^{(\mathbb{F}_2^n)^2}$ spanned all the unit vectors~$\epsilon_{\alpha,\beta}$, with~$\alpha$ and~$\beta$ words of weight~$w_1$ and~$w_2$ respectively, and~$d_H(\alpha,\beta)=d_1$.
 Then~$L_{w_1,w_2,d_1}^{(2)}$ is~$S_n$-invariant. Moreover, for any~$\lambda \vdash n$ and~$\tau \in T_{\lambda, 4}$, the irreducible representation~$\R G \cdot u_{\tau}$ is contained in~$L_{w_1,w_2, d_1}^{(2)}$, with
\begin{align}
  w_1&=  |\tau^{-1}(3)|+|\tau^{-1}(4)|, 
  \\ w_2 &=|\tau^{-1}(2)|+|\tau^{-1}(4)|, \notag 
  \\ d_1 &=|\tau^{-1}(2)|+|\tau^{-1}(3)|. \notag 
  \end{align}
  So let, for all~$\lambda \vdash n$ of height at most~$4$,
  \begin{align}
      R_{\lambda}^{(2)}:= \{\tau \in T_{\lambda, 4} \, | \,\,\,    |\tau^{-1}(2)|+|\tau^{-1}(4)|&=w, \,|\tau^{-1}(3)|+|\tau^{-1}(4)| =w,\\ |\tau^{-1}(2)|+|\tau^{-1}(3)|&\in \{0,d,d+1,\ldots,n\} \}. \notag
  \end{align}
Then 
\begin{align} \label{set2}
\{ U_{\lambda}^{(2)}:= [u_{\tau,J_2} \,\, | \,\, \tau \in  R_\lambda^{(2)} ] \,\, | \,\, \lambda \vdash n  \}
\end{align}
is representative for the action of~$S_n$ on~$\R^{(F^2)_d} \subseteq \R^{\F_2^n \times \F_2^n} $, where~$(F^2)_d$ denotes the set of all pairs of words in~$F \times F$ with distance contained in~$\{0,d,d+1,\ldots,n\}$. 

It is possible to further reduce the program (by a factor~$2$) by giving a reduction from ordered pairs to \emph{unordered} pairs of words.  We will not consider this reduction in the present paper. Regardless of a further reduction by a factor 2, the programs for computing~$A_4(n,d,w)$ are considerably larger (although they are of polynomial size) than the ones for computing~$B_4(n,d,w)$. 

Note that~$S_n$ acts trivially on~$\emptyset$. The~$S_n$-isotypical component of~$\mathbb{R}^{F^s}$ that consists of the~$S_n$-invariant elements corresponds to the matrix~$U^{(s)}_{(n)}$ in the representative set indexed by~$\lambda = (n)$. So to obtain a representative set for the action of~$S_n$ on~$(F^s)_d \cup \{\emptyset\}$ (here~$(F^1)_d := F$), we add a new unit base vector~$\epsilon_{\emptyset}$ to this matrix (as a column). 

\subsubsection{Computations for~\texorpdfstring{$D=\emptyset$}{D=empty}} 
We consider~$s=1$ and~$s=2$ for computing~$B_4(n,d,w)$ and~$A_4(n,d,w)$, respectively. If~$s=1$, then for all~$\omega \in \Omega_{2} \subseteq \Omega_4$, we define the~$\mathbb{F}_2^n \times \mathbb{F}_2^n$ matrix~$N_{\omega}^{(1)}$ by 
\begin{align}\label{blocks2}
    (N_{\omega}^{(1)})_{\alpha,\beta} := \begin{cases} 1 &\mbox{if } \{\alpha,\beta\} \in \omega,  \\ 
0 & \mbox{else.} \end{cases} 
\end{align}
Similarly, if~$s=2$, then for all 
all~$\omega \in \Omega_{4}$, we define the~$(\mathbb{F}_2^n \times \mathbb{F}_2^n) \times(\mathbb{F}_2^n \times \mathbb{F}_2^n) $ matrix~$N_{\omega}^{(2)}$ by 
\begin{align}\label{blocksundef}
    (N_{\omega}^{(2)})_{(\alpha,\beta),(\gamma,\delta)} := \begin{cases} 1 &\mbox{if } \{\alpha,\beta,\gamma,\delta\} \in \omega,  \\ 
0 & \mbox{else.} \end{cases} 
\end{align}
Let~$M_{2s,\emptyset}'(y)$ denote the matrix~$M_{2s,\emptyset}(y)$ with the row and column indexed by~$\emptyset$ removed. Then we get with~$(\ref{PhiR})$ and~$(\ref{set1})$ or~$(\ref{set2})$ that
\begin{align} \label{blocks3}
    \Phi(M_{2s,\emptyset}'(y)) = \Phi \left(\sum_{\omega \in \Omega_{2s} \setminus\{\{ \emptyset \} \}}y(\omega) N_{\omega}^{(s)} \right ) = \bigoplus_{\lambda \vdash n} \sum_{\omega \in \Omega_{2s} \setminus\{\{ \emptyset \} \}} y(\omega) {U_{\lambda}^{(s)}}^T N_{\omega}^{(s)} U_{\lambda}^{(s)}.
\end{align}
The number of~$\lambda \vdash n$ of height at most~$2^s$, and the numbers~$|R_{\lambda}^{(s)}|$, $|\Omega_{2s}|$, for~$s=1$ and~$s=2$, respectively, are all polynomially bounded in~$n$. Hence the number of blocks in~$(\ref{blocks3})$, as well as the size of each block and the number of variables occuring in all blocks are bounded by a polynomial in~$n$.
We now explain how to compute the coefficients~$u_{\tau}^T N_{\omega}^{(s)} u_{\sigma}$ in polyomial time, for~$\tau, \sigma \in R_{\lambda}^{(s)}$ and~$s \in \{1,2\}$.

If~$s=1$, define for~$P=(i_1,i_2) \in \mathbb{F}_2^{2}$, the element~$a_P := e_{i_1} \otimes e_{i_2}$, where~$e_{j}$ denotes the standard unit basis vector in~$\R^{\mathbb{F}_2}$ corresponding to~$j \in \mathbb{F}_2$. Then the set~$E^{(1)}:=\{a_P \, | \, P \in \mathbb{F}_2^{2}\}$ is a basis of~$\R^{\mathbb{F}_2} \otimes \R^{\mathbb{F}_2}$. If~$s=2$, for~$P=(i_1,i_2,i_3,i_4) \in \mathbb{F}_2^{4}$, we define $a_P := e_{i_1,i_2} \otimes e_{i_3,i_4}$,
where~$e_{i,j}$ denotes the standard unit basis vector in~$\R^{\mathbb{F}_2^2}$ corresponding to~$(i,j) \in \mathbb{F}_2^2$. Then~$E^{(2)}:=\{a_P \, | \, P \in \mathbb{F}_2^{4}\}$ is a basis of~$\R^{\mathbb{F}_2 \times \mathbb{F}_2} \otimes \R^{\mathbb{F}_2 \times \mathbb{F}_2}$. For~$s\in\{1,2\}$, let~$(E^{(s)})^*$ be the dual basis of~$E^{(s)}$.

Let~$Q^{(1)}$ and~$Q^{(2)}$ denote the sets of monomials of degree~$n$ on~$\R^{\mathbb{F}_2} \otimes \R^{\mathbb{F}_2}$ and~$\R^{\mathbb{F}_2 \times \mathbb{F}_2} \otimes \R^{\mathbb{F}_2 \times \mathbb{F}_2}$, respectively. Similar to~$(\ref{agebeuren})$, there is a natural bijection between the set~$Q^{(s)}$ and the set of~$S_n$-orbits on~$(\mathbb{F}_2^n)^{2s}$. Then the function
\begin{align}
    \left( \mathbb{F}_2^{n} \right)^{2s} &\to \mathcal{C}_{2s},\\
(x_1,\ldots,x_{2s}) &\mapsto \{x_1,\ldots,x_{2s}\}    \notag 
\end{align}
induces a surjective function $r\, : \, Q^{(s)} \to \Omega_{2s}' \setminus \{\{\emptyset\}\}$, where~$\Omega_{2s}' \supseteq \Omega_{2s}$ denotes the set of all~$S_n$-orbits of codes in~$\mathcal{C}_{2s}'$ (so not necessarily of constant weight~$w$) and~$s \in \{1,2\}$. Now define for~$\mu \in Q^{(s)}$,
\begin{align}
K_{\mu}^{(s)} =\sum_{\substack{a_1,\ldots,a_n \in E^{(s)}\\ a_1^* \cdot\ldots\cdot a_n^*=\mu} } \left( \bigotimes_{i=1}^n a_i\right) .
\end{align}
(Here every~$a_i \in (\mathbb{R}^{\F_2^s})^{\otimes 2}$ is regarded as an element of~$\mathbb{R}^{\F_2^s \times \F_2^s}$ via the natural isomorphism, so that~$K_{\mu}^{(s)}$ is an~$\mathbb{F}_2^n \times \mathbb{F}_2^n$-matrix (if~$s=1$) or an~$(\mathbb{F}_2^n \times \mathbb{F}_2^n) \times (\mathbb{F}_2^n \times \mathbb{F}_2^n)$-matrix (if~$s=2$).) By Lemma~$1$ of~$\cite{onsartikel}$ we obtain that, if~$\omega \in \Omega_{2s} \setminus\{\{ \emptyset \} \}$, then
\begin{align}
    N_{\omega}^{(s)}  = \sum_{\substack{\mu  \in Q^{(s)} \\ r(\mu) = \omega  }} K_{\mu}^{(s)}.
\end{align}
For~$\tau,\sigma \in R_{\lambda}^{(s)}$, we define
\begin{align} \label{secondp}
    p_{\tau,\sigma} := \sum_{\substack{\tau'\sim \tau_i \\ \sigma' \sim \sigma_i  }} \sum_{c, c' \in C_{\lambda}} \text{sgn}(cc') \prod_{y \in Y(\lambda)} J_s(\tau'c_i(y)) \otimes J_s(\sigma'c_i'(y)), 
\end{align}
which is a polynomial on~$\R^{\mathbb{F}_2^s} \otimes \R^{\mathbb{F}_2^s}$.
Note that each~$J_s(j) \otimes J_s(l)$ immediately gives a variable~$a_P^*$. Then with Lemma~$2$ of~$\cite{onsartikel}$, we obtain
\begin{align}
 \sum_{\mu \in Q^{(s)}} (u_{\tau}^T K_{\mu}^{(s)}u_{\sigma})\mu = p_{\tau, \sigma},
\end{align}
which is a linear combination of~$\mu\in Q^{(s)}$. Hence one can compute~$\sum_{\omega \in \Omega_{2s} \setminus\{\{ \emptyset \} \}} y(\omega)  u_{\tau}^T N_{\omega}^{(s)} u_{\sigma}$ by first expressing~$p_{\tau, \sigma}$ as a linear combination of monomials~$\mu \in Q^{(s)}$ and subsequently replacing each monomial~$\mu$ in~$p_{\tau, \sigma}$ with the variable~$y(r(\mu))$ if~$r(\mu) \in \Omega_{2s} \setminus\{\{ \emptyset \} \}$ is an orbit of minimum distance~$\geq d$ and with zero otherwise.

At last, we compute the entries in the row and column indexed by~$\emptyset$ in the matrix for~$\lambda = (n)$. Then~$\epsilon_{\emptyset}^T M_{2s,\emptyset}(y)\epsilon_{\emptyset} =M_{2s,\emptyset}(y)_{\emptyset,\emptyset}= x(\emptyset)=1$ by definition, see~$(\ref{A4ndw})$ and~$(\ref{Bndw})$. For computing the other entries we distinguish between the cases~$s=1$ and~$s=2$. If~$s=1$, then for~$\lambda = (n)$, we have~$|R^{(1)}|=1$, so there is only one coefficient to compute. If~$\tau$ is the unique $(n)$-tableau in~$\in R^{(1)}$ (containing~$w$ times symbol~$2$ and~$n-w$ times symbol~$1$), then~$u_{\tau,B} = \sum_{v \in \mathbb{F}_2^n, \text{wt}(v)=w} \epsilon_v$, so
\begin{align}\label{leeg1}
 \epsilon_{\emptyset}^T M_{2,\emptyset}(y)   u_{\tau,B} = \binom{n}{w} y(\omega_0),
\end{align}
where~$\omega_0 \in \Omega_{2}$ is the (unique) $S_n$-orbit of a code of size~$1$. 

If~$s=2$, then for~$\lambda = (n)$, any~$\tau \in R^{(2)}$ is determined by the number~$t$ of~$2$'s in the row of the Young shape~$Y((n))$ (this determines also the number of~$1$'s, $3$'s and~$4$'s). Then
\begin{align}
    u_{\tau,B} = \sum_{\substack{ v_1,v_2 \in \mathbb{F}_2^n,\\ \text{wt}(v_1)=\text{wt}(v_2)=w \\ d_H(v_1,v_2)=2t}} \epsilon_{(v_1,v_2)}. 
\end{align}
Hence
\begin{align} \label{leeg2}
 \epsilon_{\emptyset}^T M_{4,\emptyset}(y)   u_{\tau,B} = \sum_{\substack{ v_1,v_2 \in \mathbb{F}_2^n,\\ \text{wt}(v_1)=\text{wt}(v_2)=w \\ d_H(v_1,v_2)=2t}} x(\{v_1,v_2\}) = \binom{n}{w}\binom{w}{t}\binom{n-w}{t} y(\omega_t),
\end{align}
where~$\omega_t \in \Omega_{4}$ is the (unique) $S_n$-orbit of a pair of constant-weight code words of distance~$2t$.

\section{Concluding remarks}
Recently, also constant weight codes for~$n$ larger than~$28$ have been studied, see~$\cite{largen}$ and~$\cite{brouwertable}$. We therefore also provide a table with improved upper bounds for~$n$ in the range~$29 \leq n \leq 32$ and~$d \geq 10$ (cf.~Brouwer's table~$\cite{brouwertable}$). Most of the new upper bounds for these cases are instances of~$A_3(n,d,w)$, which can computed by using the block diagonalization of~$M_{2,\emptyset}(y)$ from Section~$\ref{Dempty}$ and from~$M_{3,D}(y)=M_{4,D}(y)$ for~$|D|=1$ from Section~$\ref{D2}$ (use only the blocks for~$t=0$ from this section). The bound~$A_3(n,d,w)$ is in almost all cases equal to the Schrijver bound~$\cite{schrijver}$. 

\begin{table}[ht] \footnotesize
\begin{center}
   \begin{minipage}[t]{.495\linewidth}
    \begin{tabular}{| r | r | r|| r |>{\bfseries}r | r|}
    \hline
    $n$ & $d$ & $w$ &  \multicolumn{1}{b{10mm}|}{best lower bound known}  & \multicolumn{1}{b{11mm}|}{\textbf{new upper bound}} & \multicolumn{1}{b{14mm}|}{best upper bound previously known}\\\hline 
        31& 10 & 8  & 124 & $\mathbf{322^B}$ & 329  \\    
        32& 10 & 8  & 145 & $\mathbf{402^B}$& 436 \\  
        29& 10 & 9  & 168 & $\mathbf{523^B}$& 551 \\       
        30& 10 & 9  & 203 & $\mathbf{657^B}$& 676\\   
        31& 10 & 9  & 232 & $\mathbf{822^B}$& 850 \\           
        30& 10 & 10  & 322 & 1591 &1653  \\    
        31& 10 & 10  & 465 & 2074& 2095 \\    
        32& 10 & 10  & 500 & 2669& 2720 \\    
        29& 10 & 11  & 406 & 2036& $2055^d$ \\    
        30& 10 & 11  & 504 & 2924& $2945^d$ \\    
        31& 10 & 11  & 651 & 4141& $4328^d$ \\    
        32& 10 & 11  & 992 & 5696& 6094 \\    
        29& 10 & 12  & 539 & 3091&  $3097^d$\\    
        30& 10 & 12  & 768 & 5008&  $5139^d$\\    
        31& 10 & 12  & 930 & 7259&  $7610^d$\\    
        32& 10 & 12  & 1395 & 10446& $11541^d$ \\    
        29& 10 & 13  & 756 & 4282&  $4420^d$\\    
        30& 10 & 13  & 935 & 6724&  $7149^d$\\    
        31& 10 & 13  & 1395 & 10530& $12254^d$ \\    
        32& 10 & 13  & 1984 & 16755& $18608^d$ \\    
        29& 10 & 14  & 1458 & 4927&  $5051^d$\\    
        30& 10 & 14  & 1458 & 8146&  $9471^d$\\    
        31& 10 & 14  & 1538 & 13519& $15409^d$ \\    
        32& 10 & 14  & 2325 & 22213& $24679^d$ \\       
        30& 10 & 15  & 1458 & 8948&  $10053^d$\\    
        31& 10 & 15  & 1922 & 15031& $17337^d$ \\    
        32& 10 & 15  & 2635 & 26361& $29770^d$ \\  
        32& 10 & 16  & 3038 & 27429& $30316^d$ \\          
    \hline 
    \end{tabular}\end{minipage}    
   \begin{minipage}[t]{.495\linewidth}
    \begin{tabular}{| r | r | r|| r |>{\bfseries}r | r|}
    \hline
    $n$ & $d$ & $w$ &  \multicolumn{1}{b{10mm}|}{best lower bound known}  & \multicolumn{1}{b{11mm}|}{\textbf{new upper bound}} & \multicolumn{1}{b{14mm}|}{best upper bound previously known}\\\hline 
            29& 12 & 9  &  42& $\mathbf{59^B}$& 66 \\  
        30& 12 & 9  &  42& $\mathbf{74^B}$& 94 \\ 
       31& 12 & 9  &  50& 94& 103 \\
        29& 12 & 10  &  66& 126& $129^d$ \\
        32& 12 & 11  & 186 & 573& $574^d$ \\    
        30& 12 & 12  & 190 & 492& $493^d$ \\   31& 12 & 12  & 310 & 679& $692^d$ \\     32& 12 & 12  & 496 & 952& $1014^d$ \\ 
        30& 12 & 13  & 236 & 642& $689^d$ \\ 
        31& 12 & 13  & 400 & 958& $1177^d$ \\ 
        32& 12 & 13  & 434 & 1497& $1669^d$ \\   
        29& 12 & 14 &  173& 492& $507^d$ \\ 
        30& 12 & 14 & 288 &801 & $952^d$  \\
        31& 12 & 14 & 510 &1238 & $1455^d$  \\ 32& 12 & 14 & 900 &2140 & $2143^d$  \\      
        30& 12 & 15 &  302& 894& $1008^d$ \\  
        31& 12 & 15 & 572 & 1435& $1605^d$ \\         \hline
    32& 14 & 11  & 39 & 68 & 89  \\ 
    29 & 14 & 12  & 29 & 47 & 50  \\ 
    30 & 14 & 12  & 36 & 62 & 72  \\ 
    31 & 14 & 12  & 45 & 80 & 103  \\ 
    32 & 14 & 12  & 55 & 118 & $134^d$  \\ 
    29 & 14 & 13  & 35 & 58 & 66 \\ 
    30 & 14 & 13  & 45 & 78 & $101^d$ \\ 
    31 & 14 & 13 &  60   &129    &$137^d$ \\
    29 & 14 & 14 &  58   &   63 & $82^d$ \\
    30 & 14 & 14 &  58   &   95 & $116^d$ \\
    30 & 14 & 15 &  62  &   104 & $122^d$\\    \hline  
    \end{tabular}\end{minipage}    
\end{center}
  \caption{\small An overview of the new upper bounds for constant weight codes for~$29 \leq n \leq 32$ and~$d \geq 10$. The unmarked new upper bounds are instances of~$A_3(n,d,w)$, the ones marked with~${}^B$ are instances of~$B_4(n,d,w)$. The best previously known bounds are taken from Brouwer's table~$\cite{brouwertable}$, or, when the Delsarte bound was as least as sharp, the Delsarte bound (marked with~${}^d$) is given.   \label{32}}
\end{table}

Upper bound~$(\ref{Bndw})$ could also be useful for unrestricted (non-constant weight) binary codes: one can define~$B_k(n,d)$ just as~$B_k(n,d,w)$ in~$(\ref{Bndw})$, where~$\mathcal{C}_k$ now is defined to be the collection of all codes~$C \subseteq \mathbb{F}_2^n$ of size~$\leq k$ and~$F:=\mathbb{F}_2^n$. Then one must block diagonalize~$M_{4,\emptyset}(x)$ (the block diagonalization can be found explicitly in~$\cite{semidef}$ or more conceptually in~$\cite{onsartikel}$)  and the matrices~$M_{5,D}(x)$ with~$|D|=3$ and~$|D|=2$. One can assume that~$D=\{v_1,v_2,v_3\}$ with
\begin{align}
    v_1=&\overbrace{0\ldots0\,0\ldots0}^{w}\overbrace{0\ldots0\,\,\,0\ldots0}^{n-w} \\
    v_2=& \,1\ldots1\,1\ldots1 \,0\ldots0\,\,\,0\ldots0 \notag   
    \\v_3= &\underbrace{0\ldots0}_{t_1}\underbrace{1\ldots1}_{w-t_1}\underbrace{1\ldots1}_{t_2}\underbrace{\,0\ldots0}_{n-w-t_2},    \notag 
\end{align}
for~$0 \leq w \leq n$,~$t_1\leq w$,~$t_2 \leq n-w$ and such that the weight of~$v_3$ is at least the weight of~$v_2$, so~$t_2 \geq t_1$. Then $S_{t_1} \times S_{w-t_1} \times S_{t_2} \times S_{n-w-t_2}$ acts on~$\mathcal{C}_5$, fixing~$D$ and a block diagonalization of~$M_{5,D}$ for~$|D|=3$ and~$|D|=2$ for non-constant weight codes is obtained by a straightforward adaptation of the block diagonalization of~$M_{4,D}(x)$ for~$|D|=2$ and~$|D|=1$ for constant weight codes given in Section~$\ref{D2}$.\footnote{To compute~$A_5(n,d)$, we additionally must compute a block diagonalization of~$M_{5,D}(x)$ for~$|D|=1$ (so we can assume~$D=\{0\ldots 0\}$, the zero word). This block diagonalization can be obtained by adapting the block diagonalization of~$M_{4,\emptyset}$ for constant weight codes given in Section~$\ref{Dempty}$.} However, for cases in which~$A(n,d)$ is unsettled, the program~$B_5(n,d)$ is large in practice and~$A_5(n,d)$ is still  larger (although of size polynomial in~$n$). Using a lot of computing time, one may be able to compute~$B_5(n,d)$  for some small unknown cases of~$n,d$, possibly sharpening recent semidefinite programming bounds for binary codes~$\cite{KT, semidef}$. This is material for further research.

\section*{Appendix 1: tableaux-polynomials}
We define, just as in~$\cite{onsartikel}$, for any~$n,m \in \N$,~$\lambda =(\lambda_1,\ldots, \lambda_m) \vdash  n$, and~$\tau, \sigma \in T_{\lambda,m}$, the polynomial~$p_{\tau,\sigma} \in \R[x_{i,j} \, | \, i,j=1,\ldots,m]$ by 
\begin{align}
p_{\tau,\sigma}(X)  = \sum_{\substack{\tau'\sim \tau \\ \sigma' \sim \sigma  }} \sum_{c, c' \in C_{\lambda}} \text{sgn}(cc') \prod_{y \in Y(\lambda)} x_{\tau'c(y), \sigma'c'(y)},
\end{align}
for~$X=(x_{i,j})_{i,j=1}^m \in \R^{m \times m}$. An algorithm for computing~$p_{\tau,\sigma}$ was given in~$\cite{onsartikel}$. Here we state a different method (due to Gijswijt~\cite{gijswijt}) which is very easy to implement. We define
\begin{align} 
d_{i \to j}:= \sum_{s=1}^m x_{i,s} \frac{\partial}{\partial x_{j,s}}, \,\, \text{ and }  d_{j \to i}^*:= \sum_{s=1}^m x_{s,i} \frac{\partial}{\partial x_{s,j}}.
\end{align} 
Furthermore, we write for~$i,j = 1,\ldots, m$,
\begin{align} 
t(i,j):=& \text{ $\#$ symbols~$i$ in row~$j$ of~$\tau$}, 
\\s(i,j):=& \text{ $\#$ symbols $i$ in row~$j$ of~$\sigma$}, \notag 
\end{align} 
and (here~$\lambda_{m+1}:=0$)
\begin{align} \label{pdef}
P_{\lambda}:= \prod_{k=1}^m \left( k! \, \det\left(  (x_{i,j})_{i,j=1}^{k} \right) \right)^{\lambda_{k}-\lambda_{k+1}},
\end{align}
which is a polynomial in the variables~$x_{i,j}$, where~$i,j=1,\ldots, m$. Then~$P_{\lambda}$ can be computed in time polynomially bounded in~$n$ (for fixed~$m$, note that~$\det (  (x_{i,j})_{i,j=1}^m )$ has~$m!$ terms).   
Now it holds, as is proved in~\cite[Theorem~7]{gijswijt}, that
\begin{align} \label{ultimatep}
    p_{\tau, \sigma } = \left( \prod_{j=1}^{m-1} \prod_{i=j+1}^m   \frac{1}{t(i,j)!s(i,j)!}  \right)\cdot \left( \prod_{j=1}^{m-1} \prod_{i=j+1}^m   (d_{i \to j})^{s(i,j)}  (d_{j \to i}^*)^{t(i,j)}\right) \circ P_{\lambda}.  
\end{align}
Expression~$(\ref{ultimatep})$ gives a method to compute~$p_{\tau, \sigma}$ in polynomial time (for fixed~$m$), using only methods for polynomial addition, multiplication and differentiation.

The factors~$k!$  in~$(\ref{pdef})$ can be missed for our application, i.e., for fixed~$\lambda$ one may divide the polynomial~$P_{\lambda}$, and therefore simultaneously all~$p_{\tau,\sigma}$, by $\prod_{k=1}^m ( k! )^{\lambda_{k}-\lambda_{k+1}}$ to obtain semidefinite programs with smaller numbers (but still integers) in the constraint matrices.

\section*{Appendix 2: Note regarding the computations}

In the computations for~$B_4(n,d,w)$, we added a row and a column to the block in the block diagonalization of~$M_{2,\emptyset}$ corresponding to the partition~$\lambda=(n)$. This block then gets the form
\begin{align}
T:= \begin{pmatrix}
    1 & \binom{n}{w} y(\omega_0) \\
    \binom{n}{w} y(\omega_0) & \binom{n}{w}\cdot \left(y(\omega_0) + \sum_{t=d/2}^w \binom{w}{t}\binom{n-w}{t} y(\omega_t) \right)
\end{pmatrix},
\end{align}
where the entry in the right bottom corner is computed similarly to~$(\ref{leeg1})$ and the entry in the top left corner is~$y(\omega_{\emptyset})=1$. Let~$\mathcal{M}$ be the collection of all other blocks for computing~$B_4(n,d,w)$. (So~$\mathcal{M}$ is the collection of blocks one obtains by block diagonalizing, and replacing variables~$x(C)$ by~$y(\omega)$, all matrices that are required to be positive semidefinite in~$(\ref{Bndw})$, except for block~$T$.)  Then the semidefinite program~$(\ref{Bndw})$ gets the form\footnote{Here PSD stands for \emph{positive semidefinite}.} 
\begin{align}\label{prog}
    \text{Maximize: } \binom{n}{w} y(\omega_0)\,\, \text{ subject to: $T$ PSD,~$M$ PSD for all~$M \in \mathcal{M}$} .
\end{align}
Let~$T'$ be the matrix obtained from~$T$ by dividing all entries (except for the top left entry) by~$\binom{n}{w}$. Then~$(\ref{prog})$ is equivalent to
\begin{align}\label{prog2}
    \text{Maximize: } y(\omega_0)\,\, \text{ subject to: $T'$ PSD,~$M$ PSD for all~$M \in \mathcal{M}$},
\end{align}
which is often easier for the computer to compute. Moreover, program~$(\ref{prog2})$ is equivalent to
\begin{align}\label{prog3}
    \text{Maximize: }  y(\omega_0)+ \sum_{t=d/2}^w \binom{w}{t}\binom{n-w}{t} y(\omega_t)\,\, \text{ subject to: $y(\omega_0)=1$,~$M$ PSD for all~$M \in \mathcal{M}$}.
\end{align}
To see this, write~$\text{OPT}_1$ and~$\text{OPT}_2$ for the maxima in~$(\ref{prog2})$ and~$(\ref{prog3})$, respectively. If~$(y(\omega))_{\omega \in \Omega_4}$ is a feasible assignment of the variables for~$(\ref{prog2})$, then~$(z(\omega))_{\omega \in \Omega_4}$ given by~$z(\omega) := y(\omega)/y(\omega_0)$ is a feasible assignment of the variables for~$(\ref{prog3})$ with~$z(\omega_0)+ \sum_{t=d/2}^w \binom{w}{t}\binom{n-w}{t} z(\omega_t) \geq y(\omega_0)$ (as~$T'$ is PSD). This proves~$\text{OPT}_2 \geq \text{OPT}_1$. To see the reverse inequality, suppose that~$(z(\omega))_{\omega \in \Omega_4}$ is a feasible assignment of the variables for~$(\ref{prog3})$. Then~$(y(\omega))_{\omega \in \Omega_4}$ given by
\begin{align}
 y(\omega) := z(\omega) \cdot  \left( z(\omega_0)+ \sum_{t=d/2}^w \binom{w}{t}\binom{n-w}{t} z(\omega_t) \right),
\end{align}
is feasible for~$(\ref{prog2})$ (note that~$T'$ then is positive semidefinite) with objective value $y(\omega_0)= z(\omega_0)+ \sum_{t=d/2}^w \binom{w}{t}\binom{n-w}{t} z(\omega_t)$. This proves that~$\text{OPT}_2 \leq \text{OPT}_1$ as well. Some of the new upper bounds were computed with version~$(\ref{prog3})$, and for the other upper bounds version~$(\ref{prog2})$ was used.

\section*{Appendix 3: an overview of the program}
In this section we give a high-level overview of the program, to help the reader with implementing the method. 
See Figure~$\ref{pseudocode}$ for an outline of the method.

\begin{figure}[H]
  \fbox{
    \begin{minipage}{15.4cm}
     \begin{tabular}{l}
    \textbf{Input: } Natural numbers~$n,d,w$ and~$s \in \{1,2\}$ \\
	\textbf{Output: }Semidefinite program for computing~$B_4(n,d,w)$ (if~$s=1$) and~$A_4(n,d,w)$ (if~$s=2$). \\
	$\,$\\
    \foreachv monomial~$\mu=a_{P_1}^*\ldots a_{P_n}^*$, with all~$P_i \in \F_2^4$ or all~$P_i \in \F_2^2$ (in lexicographic order)  \\
	\hphantom{1cm} assign orbit number~$r(\mu)$ to~$\mu$ (see~(\ref{orbitnumber}) and~(\ref{orbitnumber2})   below) \\
	\ndv\\
    \printv \emph{Maximize} $\binom{n}{w} y(\omega_0)$  \\
        \printv \emph{Subject to:}  \\    
\text{\color{blue}//Start with~$|D|=1$ and~$|D|=2$.}\\
\foreachv $t \in \Z_{\geq 0}$ with~$t=0$ or~$d/2 \leq t \leq w$  \\
	\hphantom{1cm} \foreachv $\bm{\lambda}=(\lambda_1,\ldots,\lambda_4) \vdash (t,w-t,t,n-w-t)$ with height$(\lambda_i) \leq 2$ for all~$i$  \\
			\hphantom{1cm} \hphantom{1cm} start a new block~$M_{\bm{\lambda}}$\\
			\hphantom{1cm} \hphantom{1cm} \foreachv $\bm{\tau} \in W_{\bm{\lambda}'}$ from~$(\ref{wlambda})$\\
					\hphantom{1cm} \hphantom{1cm}\hphantom{1cm} \foreachv $\bm{\sigma} \in W_{\bm{\lambda}}'$ from~$(\ref{wlambda})$\\
					\hphantom{1cm} \hphantom{1cm}\hphantom{1cm}\hphantom{1cm} compute~$p_{\bm{\tau},\bm{\sigma}}$ from~$(\ref{firstp})$ in variables~$a_P^*$\\					\hphantom{1cm} \hphantom{1cm}\hphantom{1cm}\hphantom{1cm}		replace monomials~$\mu$ of degree~$n$ in~$a_{P}^*$ by variables~$y(r(\mu))$ (see~(\ref{replace}) below)\\
	\hphantom{1cm} \hphantom{1cm}\hphantom{1cm}\hphantom{1cm}	 $(M_{\bm{\lambda}})_{\bm{\tau},\bm{\sigma}}:= $ the resulting linear polynomial in variables~$y(\omega)$ \\				
					\hphantom{1cm} \hphantom{1cm}\hphantom{1cm} \ndv\\
			\hphantom{1cm} \hphantom{1cm} \ndv\\
		  	\hphantom{1cm} \hphantom{1cm}   \printv $M_{\bm{\lambda}}$\emph{ positive semidefinite}.   \\
		\hphantom{1cm} \ndv 	\\   
 \ndv\\		
		\text{\color{blue}//Now~$D= \emptyset$.}\\
	\foreachv $\lambda \vdash n$ of height~$\leq 2^s$ \dov \\
		\hphantom{1cm} start a new block~$M_{\bm{\lambda}}$\\	
		\hphantom{1cm} \foreachv $\tau \in R_{\lambda}^{(s)}$ from~$(\ref{set1})$ or~$(\ref{set2})$\\
				\hphantom{1cm}\hphantom{1cm} \foreachv $\sigma \in R_{\lambda}^{(s)}$ from~$(\ref{set1})$ or~$(\ref{set2})$ \\
				\hphantom{1cm}\hphantom{1cm}\hphantom{1cm} compute~$p_{\tau,\sigma}$ from~$(\ref{secondp})$ in variables~$a_P^*$\\	\hphantom{1cm}\hphantom{1cm}\hphantom{1cm}			replace monomials~$\mu$ of degree~$n$ in~$a_{P}^*$ by variables~$y(r(\mu))$ (see~(\ref{replace}) below) \\
\hphantom{1cm}\hphantom{1cm}\hphantom{1cm}	 $(M_{\lambda})_{\tau,\sigma}:=$  the resulting linear polynomial in variables~$y(\omega)$\\				
				\hphantom{1cm}\hphantom{1cm} \ndv\\
		\hphantom{1cm} \ndv\\
		\hphantom{1cm}\ifv $\lambda=(n)$ \text{ \color{blue}// add $ \emptyset$.}\\
	\hphantom{1cm}\hphantom{1cm}	add a row and column to~$M_{\lambda}$ indexed by~$\emptyset$
	\\
	\hphantom{1cm}\hphantom{1cm}	put~$(M_{\lambda})_{\emptyset,\emptyset}:=1$ and the entries~$(M_{\lambda})_{\emptyset,\tau}$ and~$(M_{\lambda})_{\tau,\emptyset}$ as in~$(\ref{leeg1})$ or~$(\ref{leeg2})$ \\
		\hphantom{1cm}\ndv\\
		  \hphantom{1cm}\printv $M_{\lambda}$\emph{ positive semidefinite}.   \\				
	\ndv 	\\   
	\text{\color{blue}//Now~$|D|=3$ and~$|D|=4$ (i.e., nonnegativity of all variables).}\\
		\foreachv $\omega \in \Omega_4$  \\
	    \hphantom{1cm} \printv $ y(\omega) \geq 0$\\
	    \ndv 
     \end{tabular}
    \end{minipage}    } \caption{\label{pseudocode}\small{Algorithm to generate a semidefinite program for computing~$A_4(n,d,w)$ and~$B_4(n,d,w)$.}}
\end{figure}
A few remarks regarding the above steps:
\begin{enumerate}[(i)]
\item In this section we write~$\omega_t \in \Omega_{4}$ for the (unique) $S_n$-orbit of a pair of constant-weight code words of distance~$2t$, and~$\omega_{\emptyset}$ for the orbit~$\{\emptyset\}$. 
\item\label{orbitnumber} First, an orbit number~$r(\mu)$ to each monomial~$\mu = a_{P_1}^*\ldots a_{P_n}^*$, with all~$P_i \in \F_2^4$  is assigned. Each such monomial gives an~$S_n$-orbit of~$(\F_2^4)^n=(\F_2^n)^4$, and together with the map~$(\alpha,\beta,\gamma,\delta) \mapsto \{\alpha,\beta,\gamma,\delta\}$ we find a surjective function from monomials~$\mu$ to~$S_n$-orbits on~$\mathcal{C}_4'$. Each monomial that corresponds with an orbit of distance~$\geq d$ and constant weight~$w$, receives an unique orbit number~$r(\mu)$. (So that monomials~$\mu_1$ and~$\mu_2$ get the same number if and only if the monomials correspond with the same orbit.) A monomial that does not correspond with an orbit of distance~$\geq d$ and constant weight~$w$, does not get a number.
\item \label{orbitnumber2} If~$s=1$ (for computing~$B_k(n,d,w)$), then the monomials in~$p_{\tau,\sigma}$ have the form~$\mu = a_{P_1}^*\ldots a_{P_n}^*$, with all~$P_i \in \F_2^2$. So we also give these (few) monomials an orbit number.
\item \label{replace}	When replacing in~$p_{\bm{\tau},\bm{\sigma}}$ or~$p_{\tau,\sigma}$ monomials~$\mu$ of degree~$n$ in~$a_{P}^*$ by variables~$y(r(\mu))$, we only replace the monomials that got assigned a number (and hence correspond to an orbit of constant weight~$w$ and minimum distance at least~$d$). The other monomials are replaced with zero.
\end{enumerate}

\section*{Acknowledgements}
I want to thank Bart Litjens and  Lex Schrijver for very useful discussions and comments.

\selectlanguage{english} 

\end{document}